\documentclass{amsart}[12pt]
\usepackage{amsmath,amsfonts,amssymb,latexsym,amsthm}
\usepackage{graphicx,epsf}
\usepackage{epsfig}

\newtheorem{thm}{Theorem}

\newtheorem{cor}[thm]{Corollary}

\def\emp{\varnothing}

\def\nn{\mathbb N}
\def\cc{\mathbb C}

\def\Ga{\Gamma}

\def\Si{\Sigma}
\def\la{\lambda}
\def\ga{\gamma}
\def\si{\sigma}
\def\de{\delta}

\def\al{\alpha}
\def\be{\beta}
\def\om{\omega}

\def\ca{\mathcal A}

\def\CB{\mathcal B}

\def\ssu{\subset}

\def\<{\langle}
\def\>{\rangle}
\def\La{\Lambda}

\def\ts{\hskip.01cm}

\def\0{{\mathbf 0}}

\def\as{\La}

\def\sign{{\text {\rm sign} } }

\def\bi{{\mathbf{{i}}}}
\def\bi{{\mathbf{{i}}}}
\def\bj{{\mathbf{{j}}}}

\def\inv{\text{\rm inv}}

\begin{document}
\title{An algebraic extension of the MacMahon Master Theorem}

\author[Pavel~Etingof]{Pavel~Etingof$^\star$ \ }
\author[Igor~Pak]{ \ Igor~Pak$^\star$}


\thanks{\thinspace ${\hspace{-.45ex}}^\star$Department of Mathematics,
MIT, Cambridge, MA, 02139.
\hskip.06cm
Email:
\hskip.06cm
\texttt{\{etingof,pak\}@math.mit.edu}}
\maketitle

{\hskip5.83cm
July 31, 2006
}

\vskip1.3cm

\begin{abstract}
We present a new algebraic extension of the classical
MacMahon Master Theorem.  The basis of our extension
is the Koszul duality for non-quadratic algebras defined
by Berger.  Combinatorial implications are also discussed.
\end{abstract}

%

\vskip.7cm

\section*{Introduction}

The MacMahon Master Theorem is a fundamental result in
Enumerative Combinatorics, with a number of applications
to partition theory, binomial identities and q-series.
Ever since Percy MacMahon introduced it in~\cite{MM}, the
theorem was understood largely in a combinatorial and
analytic context (see~\cite{CF,GJ}).
In an important development, a quantum analogue
was discovered by Garoufalidis, L\^e and Zeilberger~\cite{GLZ},
sharply extending a basic partially commutative version
of Cartier and Foata~\cite{CF}.  Since then, several
other algebraic extensions have been introduced~\cite{FH2,KP},
but their nature remained mysterious.

Most recently, Hai and Lorenz showed that the quantum
analogue can be obtained from the Koszul duality~\cite{HL}.
While~\cite{HL} does not rederive the results of other
extensions (which appeared independently), it goes a long
way to explain them, in particular the fact that all these
algebraic extensions are based on quadratic algebras with
Koszul property.  It is a natural question whether further
(non-quadratic) extensions can be obtained.

In this paper we present an extension of the MacMahon
Master Theorem of a new type.  The construction is based
on special non-commutative algebras $\CB_{m,k}$ studied
by Berger in~\cite{B}.  These algebras have~$m$ generators
and homogenous relations of degree~$k$ coming from the exterior
algebra.  Algebras $\CB_{m,k}$ are fundamental examples of
the \emph{generalized Koszul property} for
non-quadratic algebras defined in~\cite{B}.
Our proof is uses the Hai-Lorenz approach and Berger's explicit
calculations of the Koszul resolutions in this case.
The result is a curious algebraic identity which can be stated
in a combinatorial language.

In this note, we first recall the MacMahon Master Theorem
and state the main result in Section~\ref{s:main}.  We then
give a combinatorial interpretation in Section~\ref{s:comb}.
Two natural examples and enumerative applications are
considered in sections~\ref{s:first},~\ref{s:example},
while the proofs are postponed until Section~\ref{s:proofs}.
We conclude with final remarks and open problems
(Section~\ref{s:final}).

In conclusion, let us say a few words about the fundamental
difference between an algebraic and a combinatorial
approach to the Master Theorem.  A combinatorial approach
starts by translating the algebraic identity as an equation
for two sums of words in a formal language.  This equation
is then proved by an explicit bijection.
The algebraic approach starts by showing that both
sides of the algebraic identity are invariant under a large
group of symmetries.  This action is then utilized to
obtain a concise and insightful proof of the result.
The main difference is that in the combinatorial
approach the bijection itself is of value, and
its properties often lead to new results.
We refer to~\cite{CF,FZ,KP} for examples of this
phenomenon in the context of the MacMahon Master
Theorem.  It would be interesting to obtain a
combinatorial proof of our main result (Theorem~\ref{t:main})
and further extend it in these directions.

\medskip

\section{Main results}  \label{s:main}
We begin by stating the MacMahon Master Theorem in the
classical form:

\begin{thm}  \label{t:mmt}
{\rm (MacMahon Master Theorem)} \
 Let $A=(a_{ij})_{m \times m}$ be a complex matrix,
 and let $x_1,\ldots,x_m$ be a set of variables.
 Denote by $G(k_1,\ldots,k_m)$ the coefficient of
 $x_1^{k_1}\cdots x_m^{k_m}$ in
 \begin{equation} \label{e:mmt-G}
 \prod_{i=1}^m (a_{i1}x_1+\ldots+a_{im}x_m)^{k_i}.
\end{equation}
 Let $t_1,\ldots,t_m$ be another set of variables, and
 $T = (\de_{ij} t_i)_{m \times m}$.
 Then
 \begin{equation} \label{e:mmt-id}
  \sum_{(k_1,\ldots,k_m)} \, G(k_1,\ldots,k_m)
  \ t_1^{k_1} \cdots t_m^{k_m} \, =
  \, \frac 1{\det(I - TA)},
 \end{equation}
 where the summation is over all nonnegative integer
 vectors $(k_1,\ldots,k_m)$.
\end{thm}

By taking $t_1=\ldots=t_m=1$ we get
\begin{equation} \label{e:mmt-eq}
\sum_{(k_1,\ldots,k_m)} \, G(k_1,\ldots,k_m) \, =
  \, \frac 1{\det(I - A)}\,,
\end{equation}
whenever both sides of the equation are well defined.
It is easy to see that~\eqref{e:mmt-id} and~\eqref{e:mmt-eq}
hold when all $a_{ij}$ are formal commutative variables.
Moreover, replacing $a_{ij}$ in~\eqref{e:mmt-eq} with $a_{ij}\ts t_i$
shows that~\eqref{e:mmt-eq} is  equivalent to~\eqref{e:mmt-id}, so
in the literature the MacMahon Master Theorem is often stated in
this form.

\smallskip

Let us move now into a non-commutative setting. Fix
an integer~$k$, $2\le k \le m$.
Let~$\CB =\CB_{m,k}$ be an algebra with basis
$X = \{x_1,\ldots,x_m\}$ and
the following defining relations:
\begin{equation}\label{e:CB}
\sum_{\si \in S_k} \, (-1)^{\sign(\si)} \
x_{i_{\si(1)}} \cdots \, x_{i_{\si(k)}} = \, 0,
\end{equation}
for all $1 \le i_1 < \ldots < i_k \le m$.  When $k=2$ we get
the usual symmetric algebra $\cc[x_1,\ldots,x_m]$.  Define
the set $\as(m,k)$ of \emph{admissible sequences}
$\bi= (i_1\ldots i_\ell)$, where $i_1,\ldots,i_\ell \in \{1,\ldots,m\}$,
and such that no~$k$ subsequent indices are strictly decreasing.
We call monomials
$$M_{\bi}(X) \, = \, x_{i_1}\cdots \, x_{i_\ell}\,, \ \ \
\text{where} \ \ \bi = (i_1, \ldots,i_\ell) \in \as(m,k)\,,
$$
the \emph{admissible monomials}.
It is known (see~\cite[p.~723]{B}) and easy to see that the
set of admissible monomials~$M_\bi(X)$ is a linear basis in~$\CB$
(use a Gr\"obner basis argument with the lexicographic ordering).
When $k=2$ we get the usual basis in $\cc[x_1,\ldots,x_m]$ with
indices $(i_1,\ldots,i_\ell) \in \as(m,2)$ weakly increasing.

Now fix a matrix~$A = (a_{ij})$, $a_{ij} \in \cc$, and let
$$M_\bi(A\ts X) \ = \, y_{i_1}\cdots \, y_{i_\ell}\,, \ \ \,
\text{where} \ \, y_i = a_{i1} x_1 + \ldots + a_{im} x_m, \ \,
\text{for all} \ \, 1 \le i \le m.$$
The product $M_\bi(A\ts X)$ is an element in algebra~$\CB$ and is
a generalization of the product~\eqref{e:mmt-G}.
For an admissible sequence
 $\bi=(i_1,\ldots,i_\ell) \in \as(m,k)$, denote by
 $G(\bi) = G(i_1,\ldots,i_\ell)$ the
 coefficient of~$M_{\bi}(X)$ in $M_\bi(A\ts X)$.
We set $G (\bi) = 1$ for an empty sequence $\bi = \emp$.

\begin{thm}  \label{t:main}
 Let $A=(a_{ij})_{m \times m}$ be a complex matrix,
 let $t_1,\ldots,t_m$ be another set of variables, which
 commute with $x_i$ and each other, and let
 $T = (\de_{ij} t_i)_{m \times m}$.
 Then
 \begin{equation} \label{e:main-id}
  \sum_{(i_1,\ldots,i_\ell) \in \as(m,k)} \, G(i_1,\ldots,i_\ell)
  \ t_{i_1} \cdots t_{i_\ell} \ \cdot \,
  \sum_{r \,= \,0,\,1\ \text{\rm mod}\,k} \,
  (-1)^{\al(r)} \, c_r(TA)\, = \,1\,,
 \end{equation}
 where $\al(r) = r - (r\,\ts\text{\rm mod}\,k)$,
and polynomials $c_i(TA)\in \cc[t_1,\ldots,t_m]$ are
 defined as coefficients of the {\rm characteristic polynomial:}
\begin{equation} \label{e:main-char}
\det(\la \ts I - T\ts A) \, = \, \sum_{r=0}^m \, c_r(TA) \,\la^{m-r}.
\end{equation}
\end{thm}

Note that when $k=2$ the power~$\al(r)$ is always even, and~$r$
takes all nonnegative values~$\le m$.  Therefore, the
second sum in~\eqref{e:main-id} is equal to $\det(I - TA)$,
as in~\eqref{e:mmt-id}.  Similarly, each term in the first
summation in~\eqref{e:main-id} is equal to the corresponding term
in~\eqref{e:mmt-id}.  This implies that Theorem~\ref{t:main}
is an extension of the MacMahon Master Theorem~\ref{t:mmt}.

\medskip

\section{Combinatorial reformulation} \label{s:comb}

In this section we restate our main result, Theorem~\ref{t:main},
in the language of words.  We assume that $a_{ij}$, $1\le i,j\le m$,
are commutative variables.  We begin with some combinatorial definitions.

Denote $[m] = \{1,\ldots,m\}$ and let
$\bi = (i_1,\ldots,i_\ell) \in [m]^\ell$.  Recall that
$\as(m,k)$ denotes the set of admissible sequences~$\bi$.
Denote by $\as(m,k,\ell)$ the set of sequences of length~$\ell$.
Let
$$
\inv(\bi) \, = \, \bigl|\{(s,t) : i_s > i_t, \, s < t\}\bigr|
$$
be the (usual) number of inversion in sequence~$\bi$.
Define the \emph{smallest
decreasing $k$-sequence} in~$\bi$ to be a subsequence
$i_s> i_{s+1} > \ldots > i_{s+k-1}$, where $s$ is the
smallest possible.

We say that two sequences $\bi, \ts \bi'  \in [m]^\ell$
differ by a \emph{$k$-reversion}, write $\bi \to \bi'$,
if~$\bi$ differs from~$\bi'$ by a permutation of elements
of the smallest decreasing $k$-sequence in~$\bi'$.
Of course, $\inv(\bi') > \inv(\bi)$ in this case.
For example, $(4,3,2,6,1) \to (4,6,3,2,1)$, where $k=3$,
$m=6$, $\ell=5$, and $(6,3,2)$ is the smallest decreasing
$3$-sequence in $(4,6,3,2,1) \in [6]^5$.
Here $\inv(4,3,2,6,1) = 7 < \inv(4,6,3,2,1) = 9$.

By a \emph{reversion path} from~$\bi$ to~$\bj$,
where $\bi,\bj\in [m]^\ell$,
we call a sequence of $k$-reversions:
$$\ga\,:\, \bi \, \to \, \bi' \, \to \, \bi'' \, \to \, \ldots
\, \to \, \bj
$$
We write $\bi \to_\ga \bj$ in this case, and let $|\ga|$ denotes
length (number of $k$-reversions) in path~$\ga$.
Since at every $k$-reversion the number of inversion increases,
we clearly have $\inv(\bj)-\inv(\bi) \ge |\ga|$.
Define integer coefficients
\begin{equation}\label{e:cij}
c_{\bi,\bj} \, = \, \sum_{\bi\ts\to_\ga\ts \bj} \,
(-1)^{\inv(\bj)-\inv(\bi) - |\ga|}\ ,
\end{equation}
where the summation is over all reversion paths~$\ga$
from~$\bi$ to~$\bj$.  Note that when $k=2$, we have
$c_{\bi,\bj} \in \{0,1\}$ depending on whether
$\bj \succ \bi$ in the weak Bruhat order.
Let us consider the products
$$
a_{\bi\ts\bj} \, = \, a_{i_1j_1} \,\cdots\,a_{i_\ell j_\ell},
$$
where $\bi = (i_1,\ldots,i_\ell)$,
$\bj = (j_1,\ldots,j_\ell) \in [m]^\ell$.

A \emph{partial permutation}~$\om$ of~$[m]$ is defined as a pair
$J \ssu [m]$ and the usual permutation of~$J$
(which we also denote by~$\om$).
By $\inv(\om)$ denote the (usual) number of inversion
of the permutations of elements in~$J$.
Denote by $\Si_m(r)$ the set of partial permutations of~$[m]$,
such that $|J|=r$.  For a partial permutation $\om \in \Si_m(r)$ define
$$
a_\om \, = \, a_{j_1\om(j_1)} \, \cdots \, a_{j_r\om(j_r)}\,,
$$
where $J=\{j_1,\ldots,j_r\}$.  Set $a_\om = 1$ for an empty partial
permutation (when $J = \emp$).

\smallskip


\begin{cor} \label{c:comb} For $k \le m$, in notation above
we have:
 \begin{equation} \label{e:comb-id}
  \left[1+ \sum_{\ell=1}^\infty\, \sum_{\bi \, \in \as(m,k,\ell)} \,
  \sum_{\bj \, \in [m]^\ell} \, c_{\bi,\bj} \cdot
  a_{\bi\ts\bj} \,\right] \, \cdot \,
  \left[\sum_{r \,= \,0,\,1\ \text{\rm mod}\,k} \, (-1)^{\al(r)}
  \sum_{\om \in \Si_m(r)} \, (-1)^{\inv(\om)} \, a_{\om}
  \right] \, = \, 1,
  \end{equation}
where $\al(r) = r - (r\,\ts\text{\rm mod}\,k)$.
\end{cor}

Let us emphasize that the first term in~\eqref{e:main-id} is
an infinite sum.  In contrast with the case $k=2$, it no longer
has only positive coefficients.  On the other hand, the second
term is a finite summation, and the sign pattern is completely
prescribed.

In fact, the corollary is equivalent to Theorem~\ref{t:main}.
This follows from the proof of the corollary given in
Section~\ref{s:proofs}.

\medskip

\section{First example} \label{s:first}
Let $A = (a_{ij}) = I$ be the identity matrix.
Then $M_\bi (A\ts X) = M_\bi(X)$ and $G(\bi) = 1$,
for all $\bi \in [m]^\ell$.  Therefore, the first term
in~\eqref{e:main-id} is equal to the weighted sum over all
admissible sequences, which we denote by $F_{m,k}$~:
$$F_{m,k}(t_1,\ldots,t_m) \ =
  \sum_{(i_1,\ldots,i_\ell) \in \as(m,k,\ell)} \,
 t_{i_1} \cdots t_{i_\ell}\,.
$$
For the second term, note that
$$
\det(\la\ts I - T\ts A) \, = \, (\la- t_1) \ts \cdots \ts (\la-t_m).
$$
Therefore, the coefficients of the characteristic polynomial
defined in~\eqref{e:main-char} are equal to the elementary
symmetric polynomials: $c_i = (-1)^i e_i(t_1,\ldots,t_i)$.
Now equation~\eqref{e:main-id} gives a closed formula
for the infinite summation~$F_{m,k}$~:
\begin{equation}\label{e:symm}
F_{m,k} \, = \, \frac{1}{\ts 1- e_1 + e_k - e_{k+1} + e_{2\ts k}
- e_{2\ts k\ts+1} + \, \ldots\,}\,,
\end{equation}
where indices in the denominator are~$\le m$.
In particular, this implies that $F_{m,k}$ is a symmetric
function (see Section~\ref{s:final} for further remarks).
Of course, when $k=2$, equation~\eqref{e:symm} is the classical
symmetric functions identity connecting complete and
elementary symmetric functions:
\begin{equation}\label{e:class}
1 + h_1 + h_2 + h_3 + \ts \ldots \, = \,
\frac{1}{\ts 1- e_1 + e_2 - e_{3} + \,\ldots\,+ (-1)^m \ts e_m\ts}\,.
\end{equation}
We refer to~\cite{Ma,S} for the background, proofs and other
results on symmetric functions. In the opposite extreme,
when $k=m$, we have:
\begin{equation}\label{e:case}
F_{m,m}(t_1,\ldots,t_m) \, = \,
\frac{1}{\ts 1 - (t_1+\ldots+t_m) + \ts t_1 \cdots t_m\ts}\,.
\end{equation}
This equation can be shown directly by a simple inductive
argument.

\smallskip

Consider now the enumerative applications of our equations.
Denote by \ts $L_{m,k}(\ell) = |\La(m,k,\ell)|$ \ts the
number of admissible sequences.  Taking $t_1=\ldots=t_m=t$
in~\eqref{e:symm} we obtain a generating function for
the number of admissible sequences:
\begin{equation}\label{e:genfunct}
1 \, + \, \sum_{\ell = 1}^\infty \, L_{m,k}(\ell) \, t^\ell
\, = \, \left(\ts 1 - m \ts\ts t + \binom{m}{k} \ts t^k -
\binom{m}{k+1} \ts t^{k+1} +
\ts\ldots\ts\right)^{-1}.
\end{equation}

While the exact shape of~\eqref{e:genfunct} may seem surprising,
the fact that it is rational function is not.  To see this,
consider a graph~$\Ga_{m,k}$ on $\as(m,k,m)$ with
(directed) edges of the type:
$$
(i_1,i_2,\ldots,i_m) \, \to \, (i_2,i_3,\ldots,i_m,j)\,.
$$
Clearly, $L_{m,k}(m) = m^\ell-\binom{m}{k}$, since
$(i_1 > i_2 > \ldots > i_m)$ are the
only non-admissible sequences in this case.  Observe that
admissible sequences of length~$\ell$ correspond to paths of
length $(\ell-m+1)$ in~$\Ga_m$.  Using the transfer matrix
method (see~\cite[Chapter~4]{S}) we obtain a rational
generating function for~$L_{m,k}(m)$, for any fixed~$m\ge k$.
Note also that the transfer matrix method gives denominator
of degree~$|\Ga_m|\sim m^\ell$, much larger than the degree
in~\eqref{e:genfunct}, where it is~$\le m$.

\medskip

\section{Second example} \label{s:example}
Let $k = m$ and suppose $A = (a_{ij})$, where $a_{ij} = 1$
for all $i,j\in [m]$.
In this special case we have:
$$
M_\bi (A\ts X) \, = \, (x_1 + \ldots + x_m)^\ell,
$$
for all $\bi \in [m]^\ell$.  Setting $t_1 = \ldots = t_m = t$
in~\eqref{e:main-id}, we obtain the first term:
$$
  \sum_{\bi \in \as(m,k,\ell)} \, G(\bi)
  \,\ts t^\ell \, = \, 1 \, + \, \sum_{\ell=1}^\infty \,
  N_m(\ell) \,t^\ell,
$$
where $N_m(\ell)$ is the sum of the coefficients of
all admissible monomials in $(x_1 + \ldots + x_m)^\ell$.

Now, the module condition for the second term in~\eqref{e:main-id}
are very restrictive: only $r=0$, $r=1$, and~$r=m$ are allowed.
There is a unique empty partial permutation corresponding
to $r=0$ case, and exactly~$m$ partial permutation corresponding
to $r=1$ case (all consisting of single elements).
Theorem~\ref{t:main} now gives:
$$1 \, + \, \sum_{\ell=1}^\infty \,
  N_m(\ell) \,t^\ell \, = \,
  \frac{1}{\ts |\Si_m(0)| \, - \,
  |\Si_m(1)|\,t \, + \,\det(A)\, t^m} \, =
  \, \frac{1}{\ts 1\, - \, m\ts t\ts}\,.
$$
We conclude that $N_m(\ell) = m^\ell$, which may seem quite
surprising given that the number~$L_{m,m}(\ell)$ of admissible
sequences in $\as(m,m,\ell)$ is much smaller than~$m^\ell$.
To see this, recall the generating function from the first
example:
\begin{equation}\label{e:gen-L}
1 \, + \, \sum_{\ell=1}^\infty \, L_{m,m}(\ell) \,t^\ell \,
= \,  \frac{1}{\ts 1\, - \, m\ts t \, + \,t^m\ts}\,,
\end{equation}
which is a special case of~\eqref{e:case}.  From~\eqref{e:gen-L}
it follows easily that $L_{m,m}(\ell) = \be^\ell(1+o(1))$ as
$\ell \to \infty$, for some~$\be < m$.  Of course, this can
also be seen by a direct argument.

\smallskip

Now, to explain the identity $N_m(\ell) = m^\ell$,
consider what happens when the relations~\eqref{e:CB}
in algebra~$\CB_{m,m}$ are applied to the~$m^\ell$ terms
in $(x_1 + \ldots + x_m)^\ell$.  We have:
\begin{equation} \label{e:bj}
x_{\bj} = \sum_{\bj'} \, (-1)^{\sign(\cdot)} \,x_{\bj'}\,,
\end{equation}
where the sum is over $m!-1$ permutations~$\bj'$ of the smallest
decreasing subsequence in~$\bj$ and the sign is the sign of these
permutations.  Clearly, the total sum of coefficients of terms
on the right-hand side of~\eqref{e:bj} is equal to~1.  Therefore,
the sum of coefficients of terms is unchanged under the algebra
relations, and thus equal to~$m^\ell$, the initial sum of all
coefficients.

\medskip

\section{Proof of results} \label{s:proofs}

\smallskip

\begin{proof}[Proof of Corollary~\ref{c:comb}]
First, note that we used equation~\eqref{e:main-id} with
$t_1=\ldots=t_m=1$.  The second term in~\eqref{e:comb-id}
can be easily seen to be equal to the second product
in~\eqref{e:main-id}.  This follows from the usual
expansion of determinant:
$$
\det(\la\ts I-A) \, = \ \sum_{r=0}^m \, \sum_{\om \in \Si_m(r)}
\, \la^{m-r} \,(-1)^{\inv(\om)} \,a_{\om}\,,
$$
as in~\eqref{e:main-char}.  Now, to see that the first term
in~\eqref{e:comb-id} is equal to that in~\eqref{e:main-id}
we first observe that
\begin{equation}\label{e:reversal}
x_{j_1} \, \cdots \, x_{j_\ell} \, = \,
\sum_{\bi = (i_1,\ldots,i_\ell) \in \as(m,k,\ell)} \, c_{\bi,\bj}
\ x_{i_1} \, \cdots \, x_{j_\ell}\,,
\end{equation}
for every $\bj = (j_1,\ldots,j_\ell)\in [m]^\ell$.  This follows
by induction.  Start with a product $x_{j_1} \cdots x_{j_\ell}$
corresponding to~$\bj$ and apply the relations in algebra
$\CB_{m,k}$ to the smallest decreasing $k$-sequence.  We obtain
products $x_{j_1'} \cdots x_{j_\ell'}$
corresponding to $k$-reversions $x_{\bj'} \to x_\bj$.  Now
apply this relation again, and repeat this until only admissible
sequences are obtained.  Then each reversion path $\bi \to_\ga \bj$
corresponds to a term in the summation.  After checking that the
signs in~\eqref{e:cij} correspond to those given by~\eqref{e:CB},
we obtain equation~\eqref{e:reversal}.

Note now that $G(i_1,\ldots,i_\ell)$, defined as the
 coefficient of~$M_{i_1,\ldots,i_\ell}(X)$ in $M(A\ts X)$, is equal
 to the summation of products $a_{\bi\ts\bj}$ with
 coefficients~$c_{\bi,\bj}$.  Checking the terms of the summation,
 we obtain the result.
\end{proof}

\medskip

\begin{proof}[Proof of Theorem~\ref{t:main}]
Before we move to the general case, let us first recall the
algebraic proof of the MacMahon Master Theorem (see~\cite{GLZ,HL}).
We rewrite equation~\eqref{e:mmt-eq} as follows:
$$
\sum_{\ell=0}^\infty \, {\rm Tr}\left(S^\ell A\right)
\, = \, \frac{1}{\det(I-A)}\,,
$$
where~$V$ denote the vector space of variables~$x_1,\ldots,x_m$.
Here the left-hand side is equal to the left-hand side
of~\eqref{e:mmt-eq} by definition of the symmetric power~$S^\ell(A)$.
Now, there is a natural action of $GL_m(\cc)$ on both sides.
Since matrices with distinct eigenvalues are dense in all matrices,
it suffice to show the result for diagonal matrices.  In this
case equation~\eqref{e:mmt-eq} is equivalent to
equation~\eqref{e:mmt-id} for~$A = I$.  Finally, the latter
can be written as~\eqref{e:class} which can be proved by a
straightforward calculation (see e.g.~\cite{Ma,S}).

Let us note here that $S^\ell V$ is a homogeneous component of
the natural action of $GL(V)$ on the symmetric algebra
$\ca = \cc[x_1,\ldots,x_m]$.  Then the left-hand side
of~\eqref{e:mmt-id} is exactly the character of the
$GL(V)$-action in~$\ca$.  Now Theorem~1 follows immediately
from here.

\smallskip

In the general case~$k >2$ we can proceed in a similar way
and reduce the theorem to the case of diagonal matrices,
which in turn reduces to the case~$A = I$.  As shown in
Section~\ref{s:first}, this can be written as~\eqref{e:symm}.
Unfortunately, the above mentioned straightforward
calculation is no longer possible in this case
(cf. Section~\ref{ss:id-refs}).
Following the idea in~\cite{HL}, we take a different
approach.

We use the notion of \emph{generalized Koszulity} defined
in~\cite{B}.  By Theorem~3.13 of \cite{B}, the algebras
${\mathcal B}$ defined in Section~\ref{s:main} are
generalized Koszul.  We can now proceed with the proof
of the theorem.

First, note that it is sufficient to prove the theorem in the
case when~$A$ is a nondegenerate matrix, since such matrices
are dense in the space of all matrices.  Thus, we may assume
that $A\in GL_m(\Bbb C)=GL(V)$, where $V$ (as before) is
the space of generators $x_1,\ldots,x_m$.

Since the relations of the algebra ${\mathcal B}$ are $GL(V)$-invariant,
the homogeneous components of this algebra are representations of $GL(V)$.
Therefore, the characters of these representations may be computed from
the Koszul resolution (see~\cite{Fr,PP} for the background). It follows
from~\cite{B} that
$$
{\rm Tor}_{\mathcal B}^{2j+\varepsilon}(\Bbb C,\Bbb
C)=\wedge^{kj+\varepsilon}V,
$$
sitting in degree $kj+\varepsilon$, where $\varepsilon=0,1$, and where
$\wedge^d V$ denote the exterior powers.
Therefore, the computation of the character via the Koszul resolution
yields
$$
\sum_{p\ge 0}\, {\rm Tr}\left(A|_{{\mathcal B}[p]}\right)t^p \, = \,
\left(\sum_{j\ge 0}\,\sum_{\varepsilon=0,1}\,(-1)^\varepsilon\,
{\rm Tr}\left(A|_{\wedge^{kj+\varepsilon}V}\right)\,
t^{kj+\varepsilon}\right)^{-1},
$$
where ${\mathcal B}[p]$ denotes the $p$-homogenous component.
Recall that the first factor in (\ref{e:main-id}) is exactly the
character of the $GL(V)$-action in ${\mathcal B}$.
The above equation then implies the theorem.
\end{proof}

\medskip

\section{Final remarks} \label{s:final}

\subsection{}
There are other examples of algebras with $k$-homogeneous
relations which satisfy generalized Koszul property~\cite{B}.
For example, one can omit all signs in~\eqref{e:CB} and obtain
a result similar to Theorem~\ref{t:main}, where the
second summation in~\eqref{e:main-id} has to be modified
by replacing exterior powers with symmetric powers.  By analogy
with~\eqref{e:genfunct} we get the following generating function
formula:
\begin{equation} \label{e:symm-h}
\sum_{(i_1,\ldots,i_\ell) \in \as'(m,k,\ell)} \,
 t_{i_1} \cdots t_{i_\ell}\,= \,
\frac{1}{\ts 1- h_1 + h_k - h_{k+1} + h_{2\ts k}
- h_{2\ts k\ts+1} + \, \ldots\,}\,,
\end{equation}
where $\as'(m,k,\ell)$ denote the set of sequences
$(i_1,\ldots,i_\ell) \in [m]^\ell$ without~$k$ subsequent
elements which are \emph{weakly decreasing}.  When $k=2$
we yet again get~\eqref{s:comb}, up to a substitution
$t_i \gets (-t_i)$.  For~$k>2$, both the left-hand side and
the denominator in the right-hand side are infinite series,
so the result is somewhat less natural from
combinatorial point of view.

\subsection{}
Following~\cite{GLZ} (see also~\cite{FH1,HL,KP}), there is a
natural quantum analogue of Theorem~\ref{t:main}.
The algebra relations~\eqref{e:CB} have to be replaced with
$$
\sum_{\si \in S_k} \, (-q)^{\binom{k}{2}-\inv(\si)} \
x_{i_{\si(1)}} \cdots \, x_{i_{\si(k)}} = \, 0\,,
$$
where $q \in \cc$, $q \ne 0$.  This algebra is a flat
deformation of~$\CB_{m,k}$, and its generalized Koszulity
follows from the Gr\"obner basis argument in Berger's
proof for~$\CB_{m,k}$~\cite{B}.
Another difference is in the equation~\eqref{e:main-char}
which will now contain $q$-minors.
This generalization of Theorem~\ref{t:main} is
equivariant under $U_q(GL(V))$ and can be proved along
the same lines.

One can similarly find $k$-homogenous extensions of the
Cartier-Foata Theorem and other recently considered
generalizations of the MacMahon Master Theorem
(see~\cite{FH2,KP}).

\subsection{} \label{ss:id-refs}
Although equations~\eqref{e:symm} and~\eqref{e:symm-h} do not
seem to be available in the literature, there are at least two
ways in which they can be derived from the much more general
known combinatorial results.  On one hand,
the equations follow from the \emph{Maximal String
Decomposition Theorem} (see~\cite[$\S 4.2$]{GJ}),
and on the other hand from the \emph{Cluster Theorem}
(see~\cite[$\S 2.8$]{GJ} and~\cite[Exc.~4.14]{S}).
It would be interesting to find a direct bijective
proof of the equations.

As we mentioned in Section~\ref{s:proofs}, equation~\eqref{e:symm}
and the eigenvalue argument gives an alternative proof of
Theorem~\ref{t:main}, a proof which does not use Koszulity.
Moreover, one can reverse the argument
in the proof and derive Berger's result on the Hilbert
series of $B_{m,k}$ from the theorem.  In other words, one can
now prove Berger's result combinatorially, from the Cluster or
MSD theorems.

Finally, let us mention that both Cluster and MSD theorems work in the
non-commutative setting.  In particular,  the second part
of the MSD Theorem (as stated in~\cite{GJ}) gives a natural
$q$-analogue according to the number of inversions in a
sequence.  We are curious to see the relationship between
this $q$-analogue and the quantum analogue described above.

\subsection{} To see that~$F_{m,k}$ is symmetric simply observe
that
$$
F_{m,k}(x_1,\ldots,x_m) \, = \, \sum_\nu \, s_\nu(x_1,\ldots,x_m)\,,
$$
where the summation is over all finite diagrams~$\nu$ which start
at the upper left corner and move down and to the right,
and such that all columns have size smaller than~$k$.
Such~$\nu$ are not to be confused with the rim hooks, which
are reflections of these diagrams (see e.g.~\cite{Ma,S}).
Although such~$\nu$ are not skew Young diagrams,
Schur polynomials~$s_\nu$ can be defined for any finite
arrangement of squares in~$\nn^2$.  In this case~$s_\nu$
is a summation over all sequences with the descent set
corresponding to~$\nu$.  Their symmetry follows
from the usual considerations.

\subsection{}
Even though the equation~\eqref{e:genfunct} is a simple
evaluation of~\eqref{e:symm}, the closest we could find
in the literature is the following formula:
$$1 \,+ \,
\sum_{n=1}^\infty \,e_n \, \frac{t^n}{n!} \, = \,
\left(\ts 1 \,- \,x \,+ \,\frac{x^k}{k!} \,-\,
\frac{x^{k+1}}{(k+1)!} \,+ \, \ldots\ts\right)^{-1}\,,
$$
where~$e_n$ is the number of permutations in~$S_n$
with no decreasing runs of length~$k$ (see e.g.~\cite{BD}
and~\cite[Exc. 4.2.8]{GJ}).
The connection between enumeration of permutations and
sequences with forbidden subwords is also known
(see~\cite[$\S 4.2$]{GJ} for details and further
references).

\subsection{}
As we mentioned in the introduction, the most famous
application of the MacMahon Master Theorem is to binomial
identities~\cite{GJ}.  It would be interesting to find
applications of Theorem~\ref{t:main} in this direction.

\vskip0.5cm

\subsection*{Acknowledgments}  We are grateful to Roland Berger,
Matja\v z Konvalinka, Vic Reiner and Vladimir Retakh for the
interesting discussions and help with the references.
We are especially thankful to Ira Gessel, David Jackson
and Richard Stanley who kindly showed to us how combinatorial
identities~\eqref{e:symm} and~\eqref{e:symm-h} can be derived.

\vskip1.4cm



\end{document}